\documentclass[12pt]{amsart}

\newtheorem{definition}{Definition}
\newtheorem{proposition}{Proposition}
\newtheorem{example}{Example}

\def\u{{\mathbf u}}
\def\v{{\mathbf v}}
\def\w{{\mathbf w}}

\def\<{\lbrack}
\def\>{\rbrack}

\def\operator{\mathcal}
\def\functionspace{\mathfrak}

\def\F{{\functionspace F}}
\def\P{{\functionspace P}}

\def\D{{\operator D}}
\def\A{{\operator A}}

\def\H{{\operator H}}

\def\O{{\mathcal O}}
\def\Dv{\overline{\D}v}

\def\dH{\frac{\delta\H}{\delta u}}

\def\dHone{\frac{\delta\H_1}{\delta u}}
\def\dHtwo{\frac{\delta\H_2}{\delta u}}
\def\grad{\nabla}
\def\t{^{\mathrm T}}
\def\mt{^{-\mathrm T}}
\def\m1{^{-1}}

\newcommand{\ip}[2]{\left< #1, #2 \right>}

\def\bu{{\mathbf u}}

\begin{document}

\thispagestyle{empty}
\begin{center}
{\bf Antisymmetry, pseudospectral methods, and conservative PDEs}
\end{center}
\medskip
\begin{center}
Robert I. McLachlan and Nicolas Robidoux
\end{center}
\medskip
\begin{center}
{\em Institute of Fundamental Sciences, Massey University, 
Palmerston North, New Zealand {\small \tt
R.McLachlan@massey.ac.nz, N.Robidoux@massey.ac.nz}
}
\end{center}
\bigskip
{\small\narrower\narrower
\subsection*{Abstract.} 
  ``Dual composition'', a new method of constructing energy-preserving
  discretizations of conservative PDEs, is introduced. It extends
  the summation-by-parts approach to arbitrary differential operators 
  and conserved quantities. Links to pseudospectral, Galerkin,
  antialiasing, and Hamiltonian methods are discussed.

}
\medskip

\subsection*{1. Introduction}
For all $u,v\in C^1([-1,1])$,
$$
\int_{-1}^1 v \partial_x w\, dx = 
-\int_{-1}^1 w \partial_x v\, dx + [vw]_{-1}^1,
$$
so the operator $\partial_x$ is skew-adjoint on $\{v\in C^1([-1,1]:
v(\pm1)=0\}$ with respect to the $L^2$ inner product $\ip{}{}$. Take $n$
points $x_i$, a real function $v(x)$, and
estimate $v'(x_i)$ from the values $v_i := v(x_i)$. In vector
notation, $\v' = D \v$, where $D$ is a differentiation matrix.
Suppose that 
the differentiation matrix has the form $D = S^{-1}A$, in which $S$
induces a discrete approximation
$$\ip{\v}{\w}_S := \v\t S \w\approx \int vw\,dx=\ip{v}{w},$$
 of the inner product. Then
\begin{equation}
\label{byparts}
\ip{\v}{D\w}_S + \ip{D\v}{\w}_S = \v\t S S^{-1} A \w + \v\t A\t
S\mt S \w = \v\t(A+A\t)\w,
\end{equation}
which is zero if $A$ is antisymmetric
(so that $D$ is skew-adjoint with respect to $\ip{\,}{}_S$), 
or equals $[vw]_{-1}^1$ if $x_1=-1$, $x_n=1$, and 
$A+A\t$ is zero except for $A_{nn}=-A_{11}=\frac{1}{2}$.
Eq. (\ref{byparts}) is known as a ``summation by parts'' formula;
it affects the energy flux of methods built from $D$.
More generally, preserving structural features such as skew-adjointness 
leads to natural and robust methods.

Although factorizations $D=S^{-1}A$ are ubiquitous in finite element 
methods, they have been less studied elsewhere.  They were introduced 
for finite difference methods in \cite{kr-sc} (see \cite{olsson} for 
more recent developments) and for spectral methods in \cite{ca-go}, in which
the connection between spectral collocation and Galerkin methods was used
to explain the skew-adjoint structure of some differentiation matrices.

Let $\H(u)$ be a continuum conserved quantity, the {\em energy.} 
We consider PDEs 
\begin{equation}
  \label{eq:hamilt_pde}
  \dot u = \D(u)\dH
  \mbox{,}
\end{equation}
and corresponding ``linear-gradient'' spatial discretizations
\cite{mclachlan2,mclachlan1,mqr:prl}, ODEs of
the form
\begin{equation}
  \label{eq:lin_grad}
  \dot \u = L(\u) \grad H(\u)
\end{equation}
with appropriate discretizations of $u$, $\D$, $\H$, and
$\delta/\delta u$. For a PDE of the form (\ref{eq:hamilt_pde}), if
$\D(u)$ is formally skew-adjoint, then $d\H/dt$ depends only on the
total energy flux through the boundary; if this flux 
is zero, $\H$ is an integral. Analogously, if
(\ref{eq:lin_grad}) holds, then
$\dot H = \frac{1}{2}(\nabla H)\t (L+L\t) \nabla H$,
so that $H$ cannot increase if the symmetric part of $L$ is
negative definite, and $H$ is an integral if $L$ is antisymmetric.
Conversely, all systems with an integral can be written in
``skew-gradient'' form ((\ref{eq:lin_grad}) with $L$ antisymmetric)
\cite{mqr:prl}.
Hamiltonian systems are naturally in the form
(\ref{eq:hamilt_pde}) and provide examples.

This paper summarizes \cite{mc-ro}, which contains
proofs and further examples.

\subsection*{2. Discretizing conservative PDEs}

In (\ref{eq:hamilt_pde}), we want to allow constant operators such as
$\D=\partial_x^n$ and $\D = \left(
\begin{smallmatrix}0 & 1 \\ -1 & 0\\ \end{smallmatrix}
\right)$, and nonconstant ones such as 
$\D(u) = u\partial_x + \partial_x u$.
These differ in the class of functions and boundary conditions which make 
them skew-adjoint, which suggests Defn. 1 below.

Let ($\F,\ip{}{})$ be an inner product space.  
We use two subspaces $\F_0$ and $\F_1$ which can be infinite dimensional
(in defining a PDE) or finite dimensional (in defining a discretization).
We write $\{f_j\}$ for a basis of
$\F_0$, $\{g_j\}$ for a basis of $\F_1$, and expand $u=u_j f_j$,
collecting the coefficients $(u_j)$ into a vector $\u$.
A cardinal basis is one in which $f_j(x_i) = \delta_{ij}$, so that
$u_j = u(x_j)$.

\begin{definition}
A linear operator 
$$\D: \F_0\times \F_1 \to \F,\quad \D(u)v\mapsto w\mbox{,}$$
is {\em formally skew-adjoint} if there is a functional $b(u,v,w)$,
depending only on the boundary values of $u$, $v$, and $w$ and their
derivatives up to a finite order, such that
$$
\ip{v}{\D(u)w} = -\ip{w}{\D(u)v}+b(u,v,w)\quad \forall\, u\in \F_0
,\ \forall\, v,w\in \F_1 .
$$
$\F_1$ is called a {\em domain of interior skewness} of $\D$. 
If $b(u,v,w) = 0$ $\forall\,u\in\F_0$, $\forall\,v,w\in\F_1$,
$\F_1$ is called a {\em domain of skewness} of $\D$,
and we say that $\D$ is skew-adjoint.
\end{definition}

\begin{example}\rm Let $\F^{\rm pp}(n,r) = \{u\in C^r([-1,1]):u|_{[x_i,x_{i+1}]}
\in \P_n\}$ be the piecewise polynomials of degree $n$ with $r$ derivatives.
For $\D=\partial_x$, 
$\F^{\rm pp}(n,r)$, $n,\ r\ge 0$, is a domain of interior
skewness, i.e., continuity suffices,
and $\{u\in\F^{\rm pp}(n,r):u(\pm 1)=0\}$ is a domain of skewness.
\end{example}
\begin{example}\rm
With $D(u) = 2(u\partial_x + \partial_x u) + \partial_{xxx}$, we have
$$
\ip{v}{\D(u)w}+\ip{w}{\D(u)v} = [w_{xx}v - w_x v_x + w v_{xx} + 2 uvw],$$
so suitable domains of interior skewness are $\F_0 = \F^{\rm
  pp}(1,0)$, $\F_1=\F^{\rm pp}(3,2)$, i.e., more smoothness is required
from $v$ and $w$ than from $u$.
A boundary condition which makes $\D(u)$ skew is $\{v:
v(\pm 1)=0,\ v_x(1)=v_x(-1) \}$.
\end{example}

\begin{definition} $\F_0$ is
  {\em natural for $\H$} if $\forall u \in \F_0$ there exists
  $\frac{\delta \H}{\delta u}\in\F$ such that
\[
 \lim_{\varepsilon\rightarrow 0} 
 \frac{ \H(u+\varepsilon v) - \H(u) }{ \varepsilon }
 = \ip{v}{\frac{\delta \H}{\delta u}}
 \quad \forall\, v\in\F
 \mbox{.}
\]
\end{definition}

The naturality of $\F_0$ often follows from the vanishing of the
boundary terms, if any, which appear of the first variation of $\H$,
together with mild smoothness assumptions.

We use appropriate 
spaces $\F_0$ and $\F_1$ to generate spectral, pseudospectral, and
finite element discretizations which have discrete energy
$H:=\H|_{\F_0}$ as a conserved quantity.  The discretization of the
differential operator $\D$ is a linear operator $\overline\D
:\F_1\to\F_0$, and the discretization of the variational derivative
$\dH$ is $\overline{\dH}\in\F_1$. 
Each of $\overline \D$ and $\overline\dH$ is a weighted residual
approximation \cite{finlayson}, but each uses spaces of 
weight functions different from its space of trial functions.  

\begin{definition}
$S$ is the matrix of $\ip{}{}|_{\F_0\times\F_1}$, i.e. 
$S_{ij} := \ip{f_i}{g_j}$.
$A(u)$ is the matrix of the linear operator $\A:(v,w)\mapsto\ip{v}{\D(u)w}$,
i.e. $A_{ij}(u) := \ip{g_i}{\D(u)g_j}$.
\end{definition}

\begin{proposition}
  Let $\F_0$ be natural for $\H$ and let $S$ be nonsingular. Then for
  every $u\in\F_0$ there is a unique element
  $\overline{\frac{\delta \H}{\delta u}}\in\F_1$ such that
\[
\ip{w}{\overline{\frac{\delta \H}{\delta u}}} =
\ip{w}{\frac{\delta \H}{\delta u}} \quad \forall\, w\in\F_0
\mbox{.}
\]
Its coordinate representation is $S\m1\nabla H$ where $H(\u):=\H(u_i f_i)$.
\end{proposition}

\begin{proposition}
  \label{prop:D}
  Let $S$ be nonsingular.  For every $v\in\F_1$, there exists a 
  unique element $\Dv\in\F_0$ satisfying
    \[
    \ip{\Dv}{w} = \ip{\D v}{w} \quad \forall\, w\in\F_1 \mbox{.}
    \]
  The map $v\mapsto\Dv$ is linear, with matrix representation $D:=S\mt A$.
\end{proposition}

\begin{definition}
  $\overline{\D}\overline{\frac{\delta\H}{\delta u}}:\F_0\to\F_0$
  is the {\em dual composition discretization} of
  $\D\frac{\delta\H}{\delta u}$.
\end{definition}

Its matrix representation is $S\mt A S\m1 \nabla H$.
The name ``dual composition'' comes from the dual roles played
by $\F_0$ and $\F_1$ in defining $\overline{\D}$ 
and $\overline{\dH}$
which is necessary so that their composition has the required
linear-gradient structure.
Implementation and accuracy of
dual composition and Galerkin discretizations are similar. Because
they coincide in simple cases, such methods are widely used already.

\begin{proposition}
  If $\F_1$ is a domain of skewness, the matrix $S\mt A S\m1$
  is antisymmetric, and the system of ODEs
  \begin{equation}
  \label{eq:disc}
  \dot\bu
  =
  S\mt A S\m1 \nabla H
  \end{equation}
  has $H$ as an integral. If, in addition, $\D$ is constant---i.e.,
  does not depend on $u$---then the system (\ref{eq:disc}) is Hamiltonian.
\end{proposition}

The method of dual compositions also yields 
discretizations of linear differential operators $\D$ (by taking
$\H=\frac{1}{2}\ip{u}{u}$), and discretizations of variational
derivatives (by taking $\D=1$).
It also applies to formally {\em self}-adjoint
$\D$'s and to mixed (e.g. advection-diffusion) operators, where
preserving symmetry gives control of the energy.

The composition of two weighted residual discretizations is not
necessarily itself of weighted residual type. The simplest case is
when $\F_0=\F_1$ and we compare the dual composition to the 
{\em Galerkin discretization}, a weighted
residual discretization of $\D \frac{\delta \H}{\delta u}$ with
trial functions and weights both in $\F_0$. They are the same when
projecting $\dH$ to $\F_0$, applying $\D$, and
again projecting to $\F_0$, is equivalent to directly projecting
$\D\dH$ to $\F_0$.

For brevity, we assume $\F_0=\F_1$ for the rest of Section 2.

\begin{proposition}
  \label{prop:galerkin}
$\overline{\D}\overline{\dH}
    $ is the Galerkin approximation of
  $\D \dH$ if and only if
$ \D \big( \overline{\dH} - \dH \big) \perp \F_0.$
This occurs if 
  (i) $\D(\F_0^\perp)\perp\F_0$, or
  (ii) $\overline\D$ is exact and applying $\D$ and orthogonal
	projection to $\F_0$ commute, or
  (iii) $\overline{\frac{\delta \H}{\delta u}}$ is exact,
i.e., $\dH\in\F_0$.
\end{proposition}

Fourier spectral methods with $\D=\partial_x^n$ satisfy (ii), since 
then $\F$ has an orthogonal
basis of eigenfunctions ${\mathrm e}^{ijx}$ of $\D$, and differentiating
and projecting (dropping the high modes) commute. This is illustrated
later for the KdV equation.

The most obvious situation in which $\dH\in\F_0$ is when
$\H=\frac{1}{2}\ip{u}{u}$, since then
$\dH=u\in\F_0$ and $\D\frac{\delta \H}{\delta u}=\D u$,
and the discretization of $\D$ is obviously the Galerkin one! 
When the functions $f_j$ are nonlocal, $D$ is often called the
spectral differentiation matrix.  The link to standard pseudospectral
methods is that some Galerkin methods are pseudospectral.

\begin{proposition}
\label{prop:pseudo}
If $\D(\F_1)\subseteq\F_1$, then $\overline{\D}v=\D v$, 
i.e., the Galerkin approximation of the derivative is exact.
If, further, $\{f_j\}$ is a cardinal basis,
  then $D$ is the standard pseudospectral differentiation matrix,
i.e. $D_{ij} = \D f_j(x_i)$.
\end{proposition}

We want to emphasize that although $A$, $S$, and $D$ depend on the basis, 
$\overline\D$ depends only on $\F_0$ and $\F_1$, i.e., it is
basis and grid independent. 
In the factorization $D=S\mt A$, the (anti)symmetry of $A$ and $S$ is basis
independent, unlike that of $D$. These points are well known in 
finite elements, less so in pseudospectral methods.

\begin{example}[\bf Fourier differentiation\rm]\rm
Let $\F_1$ be the trigonometric polynomials of degree $n$, which is
closed under differentiation (so that Prop. \ref{prop:pseudo}) applies,
and is a domain of skewness of $\D=\partial_x$. In any basis, $A$ is
antisymmetric. Furthermore, the  two popular bases, $\{\sin(j x)_{j=1}^n,
\cos(j x)_{j=0}^n\}$, and the cardinal basis on equally-spaced grid
points, are both orthogonal, so that $S=\alpha I$ and $D=S^{-1}A$ is
antisymmetric in both cases.
\end{example}

\begin{example}[\bf Polynomial differentiation\rm]\rm
\label{sec:cheb}
$\F_1=\P_n([-1,1])$ is a domain of interior skewness which is
closed under $\D=\partial_x$, so pseudospectral differentiation
factors as $D=S^{-1}A$ in any basis. For a cardinal
basis which includes $x_0=-1$, $x_n=1$, we have $(A+A\t)_{ij}=-1$
for $i=j=0$, $1$ for $i=j=n$, and 0 otherwise, making obvious
the influence of the boundary.
For the Chebyshev points $x_i = -\cos(i
\pi/n)$, $i=0,\dots,n$, $A$ can be evaluated first in a basis
$\left\{ T_i \right\}$ of Chebyshev polynomials:
one finds
$A_{ij}^{\rm cheb} = 2 j^2/(j^2-i^2)$ for $i-j$ odd, and
$S_{ij}^{\rm cheb} -2(i^2+j^2-1)/
    [((i+j)^2-1)((i-j)^2-1)]$ for $i-j$ even, with other entries 0.
Changing to a cardinal basis by
$F_{ij} = T_j(x_i) = \cos(i j \pi/n)$, a
discrete cosine transform, gives $A=F^{-1} A^{\rm cheb} F\mt$.
For example, with $n=3$ 
(so that $(x_0,x_1,x_2,x_3)=(-1,-\frac{1}{2}, \frac{1}{2},1)$), we have
$$ D = 
  {\scriptstyle \frac{1}{6}}
  \left(
    \begin{smallmatrix}
      -19 & 24 & -8 & 3 \\
      2 & -6 & -2 & 6 \\
      -6 & 2 & 6 & -2 \\
      -3 & 8 & -24 & 19 \\
    \end{smallmatrix}
  \right)
  = S\mt A = 
  {\scriptstyle \frac{1}{256}}
  \left(
    \begin{smallmatrix}
      4096 & -304 & 496 & -1024\\
      -304 & 811 & -259 & 496\\
      496 & -259 & 811 & -304\\
      -1024 & 496 & -304 & 4096\\
    \end{smallmatrix}
  \right)
  {\scriptstyle \frac{1}{270}}
  \left(
    \begin{smallmatrix}
      -135 & 184 & -72 & 23 \\
      -184 & 0 & 256 & -72\\
      72 & -256 & 0 & 184 \\
      -23 & 72 & -184 & 135 
    \end{smallmatrix}
  \right).
$$
$S$ and $A$ may be more amenable to study than $D$ itself. 
All their eigenvalues are very well-behaved; none are spurious. The
eigenvalues of $A$ are all imaginary and, as $n\to\infty$, uniformly
fill $[-i\pi,i\pi]$ (with a single zero eigenvalue corresponding
to the Casimir of $\partial_x$). 
The eigenvalues of $S$ closely approximate the
quadrature weights of the Chebyshev grid.
\end{example}

For $\D\ne\partial_x$, $\overline{\D}$ may be quite expensive
and no longer pseudospectral. (There is in general no
$S$ with respect to which the pseudospectral approximation of
$\D v$ is skew-adjoint.) However, $\overline{\D}v$ can 
be computed quickly if fast transforms between cardinal and
orthonormal bases exist. We evaluate $\D v$ exactly for
$v\in\F_1$ and then project $S$-orthogonally to $\F_1$.

\begin{example}[\bf Fast Fourier Galerkin method\rm]\rm
\label{fastfourier}
  Let $\D(u)$ be linear in $u$, for example, $\D(u) = u\partial_x
  + \partial_x u$. Let $u,\ v\in\F_1$, the trigonometric
  polynomials of degree $n$. Then $\D(u)v$ is
  a trigonometric polynomial of degree $2n$, the first $n$ modes of
  which can be evaluated exactly using antialiasing and Fourier
  pseudospectral differentiation. The approximation whose error
  is orthogonal to $\F_1$ is just these first $n$ modes, because $S=I$
  in the spectral basis. That is, the antialiased
  pseudospectral method is here identical to the Galerkin method, and hence
  skew-adjoint. Antialiasing makes pseudospectral methods conservative.
  This is the case of the linear $\D$'s of the Euler fluid equations.
\end{example}

\begin{example}[\bf Fast Chebyshev Galerkin method\rm]\rm
  Let $\D(u)$ be linear in $u$ and let $u,\ v\in\F_1=\P_n$.
  With respect to the cardinal basis on the Chebyshev grid with $n+1$
  points, $\overline{\D}(u)v$ can be computed in time $\O(n \log n)$ as follows:
(i)
   Using an FFT, express $u$ and $v$ as Chebyshev polynomial
    series of degree $n$;
(ii)   Pad with zeros to get Chebyshev polynomial series of formal
    degree $2n$;
(iii)   Transform back to a Chebyshev grid with $2n+1$ points;
(iv)   Compute the pseudospectral approximation of $\D(u)v$ on the
    denser grid. Being a polynomial of degree $\le 2n$, the
    corresponding Chebyshev polynomial series is exact;
(v)   Convert $\D(u)v$ to a Legendre polynomial series using a fast
    transform \cite{al-ro};
(vi)   Take the first $n+1$ terms. This produces 
$\overline{\D}(u)v$, because the Legendre
polynomials are orthogonal.
(vii)   Convert to a Chebyshev polynomial series with $n+1$ terms
    using a fast transform;
(viii)   Evaluate at the points of the original Chebyshev grid using an FFT.
\end{example}

\subsection*{3. Examples of the dual composition method}
\begin{example}[\bf The KdV equation\rm]\rm
$ \dot u + 6 u u_x + u_{xxx}=0$ with
periodic boundary conditions has features which can be used to illustrate 
various properties of the dual composition method. Consider two of its
Hamiltonian forms,
$$
\dot u = \D_1\dHone\mbox{, } \D_1 =
\partial_x\mbox{, } \H_1 = \int\big( -u^3+\frac{1}{2} u_x^2\big)\,dx\mbox{,}$$
and
$$
\dot u = \D_2\frac{\delta \H_2}{\delta u}\mbox{, } \D_2 =
-(2u\partial_x + 2\partial_x u + \partial_{xxx})\mbox{, } \H_2 =
\frac{1}{2}\int u^2\,dx\mbox{.}$$

In the case $\F_0=\F_1=\F^{\rm trig}$, $v:=\overline{\dHone}$ 
is the orthogonal projection to $\F_0$ of $\dHone=-3u^2-u_{xx}$; this can be
computed by multiplying out the Fourier series and dropping all but
the first $n$ modes, or by antialiasing.
Then $\overline{\D}_1 v = v_x$, since
differentiation is exact in $\F^{\rm trig}$.  Since $\D_1$
is constant, the discretization is a Hamiltonian system, and since
$\overline{\D}_1$ is exact on constants, it also
preserves the Casimir ${\operator C}=\int u\,dx$.
In this formulation, Prop. \ref{prop:galerkin} (ii) shows that
the dual composition and Galerkin approximations of $\D_1\dHone$ coincide,
for differentiation does not map high modes to lower modes, i.e., 
$\D_1(\F^{{\rm trig}\perp})\perp\F^{\rm trig}$.

In the second Hamiltonian form, $H_2 = \frac{1}{2}\u\t S \u$, $\dHtwo = 
S^{-1}\nabla H_2 = \u,$ and the Galerkin approximation of $\dHtwo$ is exact, 
so that Prop. \ref{prop:galerkin} (iii) implies that the composition 
$\overline\D_2\overline{\dHtwo}$ {\em also} coincides with the Galerkin 
approximation.  $\overline{\D}_2v$ can evaluated using antialiasing 
as in Example \ref{fastfourier}.  $\overline{\D}_2$ is
not a Hamiltonian operator, but still generates a skew-gradient
system with integral $H_2$. Thus in this (unusual) case,
the Galerkin and antialiased pseudospectral methods coincide and have
three conserved quantities,
$H_1$, $H_2$, and ${\operator C}|_{\F^{\rm trig}}$.

The situation for finite element methods with 
$\F_0=\F_1=\F^{\rm pp}(n,r)$ is different.
In the first form, we need $r\ge 1$ to ensure
that $\F_0$ is natural for $\H_1$; in the second form, naturality is
no restriction, but we need $r\ge2$ to ensure that $\F_1$ is a domain
of interior skewness.  The first dual composition method 
is still Hamiltonian with
integral $H_1$ and Casimir $C=u_i\int f_i\, dx$, but because
$\overline\D_1$ does not commute with projection to $\F_1$, it is {\em
not} a standard Galerkin method. 
In the second form, $\dHtwo=u$ is still exact, so the
dual composition and Galerkin methods still coincide.
However, they are not Hamiltonian.
\end{example}

\begin{example}[\bf An inhomogeneous wave equation\rm]\rm
When natural and skew boundary conditions conflict, it is necessary
to take $\F_0\ne\F_1$. Consider
$ \dot q  = a(x)p$, $\dot p = q_{xx}$, $q_x(\pm1,t)=0$.
This is a canonical Hamiltonian system with
$$ \D = \left(\begin{matrix}0 & 1 \\ -1 & 0 \\\end{matrix}\right),\ 
 \H = \frac{1}{2}\int_{-1}^1 \big(a(x)p^2 + q_x^2\big)\, dx,\ 
\frac{\delta \H}{\delta q} = -q_{xx},\ 
\frac{\delta \H}{\delta p} = a(x)p.$$
Note that (i) the boundary condition is
natural for $\H$, and (ii) 
no boundary conditions are required for $\D$ to be skew-adjoint in $L^2$.
Since $\overline{\dH}$ is computed with trial functions in $\F_1$, we
should not include $q_x(\pm1)=0$ in $\F_1$, for this would be to 
enforce $(-q_{xx})_x=0$.
In \cite{mc-ro} we show that a spectrally accurate dual composition method is 
obtained with
$ \F_0 = \{ q\in \P_{n+2}: q_x(\pm 1)=0 \} \times \P_n$ and
$  \F_1  = \P_n\times \P_n$.
\end{example}

\subsection*{4. Quadrature of Hamiltonians}
\label{sec:quadrature}
Computing $\nabla H =\nabla\H(u_j f_j)$ is not always possible in closed form. 
We would like to approximate $\H$ itself by quadratures in real space.
However, even if the discrete $H$ and its gradient are spectrally accurate
  approximations, they cannot always be used to construct spectrally
  accurate Hamiltonian discretizations.

In a cardinal basis, 
let $\H=\int h(u)dx$ and define the
quadrature Hamiltonian $H_q:= h( u_j) w_j = \w\t h(\u)$
where $w_j = \int f_j dx$ are the quadrature weights.
Since $\nabla H_q = W h'(\u)$, $\dH\approx W^{-1}\nabla H_q$, 
Unfortunately,
$DW^{-1}\nabla H_q$ is not a skew-gradient system, while
$D S^{-1} \nabla H_q$ is skew-gradient, but is not an accurate approximation.

$D W^{-1} \nabla H_q$ can only be a skew-gradient
system if $DW^{-1}$ is antisymmetric, which occurs in three general cases.
(i) On a constant grid, $W$ is a multiple of the identity, so
  if $D$ is antisymmetric, $D W^{-1}$ is too.
(ii) On an arbitrary grid with $D=\left(
\begin{smallmatrix}
  0 & I \\
  -I & 0\\
\end{smallmatrix}\right)$,
$DW^{-1}$ is antisymmetric.
(iii) On a Legendre grid with $\F_0=\F_1$,
 $S=W$, and $D W^{-1} = W^{-1} A W^{-1}$ is antisymmetric.
The required compatibility between
$D$ and $W$ remains an intriguing and frustrating obstacle to the
systematic construction of conservative discretizations of strongly
nonlinear PDEs.


\end{document}